\documentstyle[12pt,subeqnarray,verbatim]{article}

\begin{document}
\newcommand{\dis}{\displaystyle}
\newcommand{\dml}{\rm dim ~ }
\newcommand{\krl}{\rm ker ~ ~ }
\newcommand{\expon}{\rm e}
\newcommand{\id}{ 1 \hspace{-2.85pt} {\rm I} \hspace{2.5mm}}
\newcommand{\func}{ (\frac{\Phi(\Gamma) - \Phi(J_0)}{\Gamma -
[J_0][J_0 + 1]})^{\frac{1}{2}} }
\newcommand{\Dlt}{\Delta}
\newcommand{\nJ}{\widehat{J}}
\newcommand{\ca}{{\cal C}}
\newcommand{\ct}{\tilde{\cal C}}
\newcommand{\Jp}{\tilde{J}_+}
\newcommand{\Jm}{\tilde{J}_-}
\newcommand{\Jz}{\tilde{J}_0}
\newcommand{\Jpm}{\tilde{J}_{\pm}}
\newcommand{\U}{{\cal U}}
\newcommand{\alp}{\alpha^\prime}
\newcommand{\be}{\begin{equation}}
\newcommand{\ee}{\end{equation}}
{\thispagestyle{empty}
\rightline{} %Paper Code Number
\rightline{} %Code Number
\rightline{} %Month, year
\vskip 1cm
%%%%%%%%%%%%%%%%%%%%%%%%  title %%%%%%%%%%%%%%%%%%%%%%%%%%%
\centerline{\large \bf Some comments on }
\centerline{\large \bf $q$-deformed oscillators}
\centerline{\large \bf and $q$-deformed $su(2)$ algebras}

%%%%%%%%%%%%%%%%%%%%%%%%%%%%%%%%%%%%%%%%%%%%%%%%%%%%%%%%%%%

\vskip 2cm
%%%%%%%%%%%%% Author  and Address %%%%%%%%%%%%%%%%%%%%%%%%%
\centerline{L. C. Kwek 
{\footnote{E-mail address:
scip6051@leonis.nus.edu.sg }},
and C. H. Oh 
{\footnote{E-mail address:
phyohch@leonis.nus.edu.sg }}
}
\centerline{{\it Department of Physics, Faculty of Science, } }
\centerline{{\it National University of Singapore, Lower Kent Ridge,} }
\centerline{{\it Singapore 119260, Republic of Singapore. } }}%address
\vskip 0.1in
%%%%%%%%%%%%%%%%%%%%%%%%%%%%%%%%%%%%%%%%%%%%%%%%%%%%%%%%%%%

\vskip 1cm
\centerline{\bf Abstract} \vspace{5mm}

\noindent{The various relations between $q$-deformed
oscillators algebras  and the $q$-deformed $su(2)$
algebras are discussed.  
In particular, we exhibit 
the similarity of the
$q$-deformed $su(2)$ algebra obtained from $q$-oscillators via
Schwinger construction and those obtained from $q$-Holstein-Primakoff
transformation and show how the relation between $su_{\sqrt{q}}(2)$ and
Hong Yan $q$-oscillator can be regarded
as an special case of In\"{o}u\"{e}- Wigner contraction.
This latter observation and the imposition of positive norm requirement
suggest that Hong Yan $q$-oscillator algebra is different from the usual
$su_{\sqrt{q}}(2)$ algebra, contrary to current belief in the literature.

\newpage
}
%%%%%%%%%%%%%%%%%%%%%%% Introduction %%%%%%%%%%%%%%%%%%%%%%%%
\section{Introduction}
Since the Macfarlane-Biedenharn (MB) papers \cite{mac,bied}
on the construction of $su_q(2)$
algebra from the q-deformed oscillator algebra \`a la Schwinger way,
there are by now many different versions of the q-deformed algebra. 
However all these $q$-deformed oscillator algebras are not Hopf algebras
except the Hong Yan type and its generalization \cite{yan,oh}.
It should be stressed
here that via the Schwinger construction, it is only the `algebraic'
aspect of the Hopf algebra $su_q(2)$ which can be expressed  in terms of
the $q$-oscillator algebra; the co-algebraic structure of $su_q(2)$
cannot be easily obtained from the q-oscillator algebra granted that the
latter possesses a Hopf structure.

It has been claimed \cite{cel1,oh2} that Hong Yan (HY) Hopf algebra 
is the same as the
$su_q(2)$ Hopf algebra and a formal relation has been established for
the generators of $su_{\sqrt{q}}(2)$ and the HY oscillator algebra.
Nevertheless if we impose positive norm requirement for the states, then
at the representation level, the identification breaks down for some
values of $|q| =1$, since for these values, the positive norm
requirement does not hold.  In fact, the positive norm requirement
\cite{fuji1} is in
conflict with the truncation condition \cite{oh2}
imposed on the states of the
oscillator so as to get finite multiplets for $su_{\sqrt{q}}(2)$.  In
other words, for $|q| =1$ ($q= \expon^{i \epsilon}, ~ ~ \epsilon$
arbitrary) HY oscillator algebra is different from
$su_{\sqrt{q}}(2)$ algebra.  Furthermore, although $su_{\sqrt{q}}(2)$
has a $q \rightarrow 1$ limit at the coalgebra level, the coalgebraic
structure for HY fails in this limit. In the following section, we
summarize the $q$-Schwinger construction of $q$-deformed $su(2)$ algebra
in terms of a pair of $q$-oscillator algebras; different $q$-oscillator
algebras lead to different $q$-deformed $su(2)$ algebras.  Most authors
prefer to set the Casimir in their $q$-Schwinger construction to zero.
However, one sometimes find it convenient and essential to consider
{\it non-zero} Casimir for some physical applications\cite{fuji1,lorek}.
A natural generalization with two additonal parameters $\alpha$ and
$\beta$ is also provided.
In section \ref{qhols}, 
we exhibit results for $q$-Holstein Primakoff (HP) transformation with
{\it non-zero} Casimirs for the MB and HY oscillators. 
The results are similar to those presented in section \ref{qcontr}.
Different contractions of $q$-deformed $su(2)$ algebras to
the various $q$-oscillator algebras are elucidated in section 
\ref{sect4}.  In particular, we show that the relation between
$su_{\sqrt{q}}(2)$ and HY $q$-oscillator algebras obtained in 
ref \cite{cel1,oh2} can be regarded as a form of contraction. 
In the last section, we point out explicitly that
at the representation level
the usual
$su_{\sqrt{q}}(2)$ algebra is not the same as the HY
$q$-oscillator algebra.

We recall that the
quantum universal enveloping algebra,
${\cal U}_q(su(2))$, was first studied by Skylanin \cite{sky} and
independently by Kulish and Reshetikhin \cite{kul}. 
This algebra has been applied extensively to the study of
the eight vertex models,
the $XXZ$ ferromagnetic and anti-ferromagnetic models and the
sine-Gordon models. The universal enveloping algebra, ${\cal U}_q(su(2))$
is generated by three operators, $J_{\pm}$ and $J_0$ satisfying the
commutation relations
\begin{subeqnarray}
\lbrack J_0, J_{\pm} \rbrack & = & \pm J_{\pm}, \\
\lbrack J_+, J_- \rbrack & = & \lbrack 2J_0 \rbrack, \label{cr}
\end{subeqnarray}
where $[x]$ denotes $\dis \frac{q^x - q^{-x}}{q - q^{-1}}$.

A generalized $q$-deformed $su(2)$ algebra \cite{bonat,poly}
has also been proposed 
and the operators
$\nJ_{\pm}$ and $\nJ_0$ satisfies a modified 
commutation relations
\begin{subeqnarray}
\lbrack \nJ_0, \nJ_{\pm} \rbrack & = & \pm \nJ_{\pm}, \\
\lbrack \nJ_+, \nJ_- \rbrack & = & \Phi(\nJ_0(\nJ_0 + 1)) -
\Phi(\nJ_0(\nJ_0 - 1)), \\
& = & \Psi(\nJ_0) -
\Psi(\nJ_0 - 1), \label{mcr}
\end{subeqnarray}
where the functions $\Phi(\nJ_0)$ and $\Psi(\nJ_o)$ are some suitably
chosen functions of $\nJ_0$. It has been shown in ref\cite{poly} that
the imposition of hermiticity condition requires the generalized $q$-deformed
$su(2)$ to assume the form given in eq(\ref{mcr}).

\section{$q$-Schwinger Construction} \label{qcontr}
Traditionally, the algebra $su(2)$ can be realized in terms of a pair of
bosonic
creation and annihilation operators of a harmonic oscillator using the
Schwinger construction.  A $q$-analogue of this construction is given by
MB\cite{mac,bied,kulish}.
The operators $a, a^\dagger$ and $N$ of the $q$-deformed
oscillator
algebra obey the relations
\begin{subeqnarray}
\lbrack N, a^\dagger \rbrack =  a^\dagger, & &
\lbrack N, a \rbrack = - a, \\
a a^{\dagger} - q^{-1} a^{\dagger} a & = & q^{N}, \\
\ca_1 & = & q^{N} (a^\dagger a - [N]). \label{macf}
\end{subeqnarray}
This oscillator algebra does not appear to possess
a Hopf structure. But a Hopf structure is possible
for another version of  the $q$-deformed oscillator which was first
proposed by HY \cite{yan} in which the operators $a, a^\dagger$ and
$N$ satisfy eq(\ref{macf}a) and eq(\ref{macf}b) and
\begin{subeqnarray}
\lbrack a, a^{\dagger} \rbrack & = & [N + 1] - [N] \\
\ca_{2} & = & a^\dagger a - [N]. \label{hyo}
\end{subeqnarray}
In general, these two versions of the $q$-deformed oscillator algebras
are not equivalent\cite{oh} although the
two algebras coincide on the usual `Fock' space basis $|n>$.  

Mathematically, it has always been intrinsically appealing and insightful 
to generalize a particular mathematical structure as much as 
possible\cite{duc,rideau,quesne}. 
One possible generalization of 
the MB algebra is to introduce two 
additional parameters $\alpha$ and $\beta$.
One then defines the generalized MB (GMB) algebra\cite{duc}
with the relations
in eq(\ref{macf}b) and eq(\ref{macf}c) replaced by
\begin{subeqnarray}
a a^{\dagger} - q^{\alpha} a^{\dagger} a & = & q^{\beta N}, \\
\ca_3 & = & q^{-\alpha N} (a^\dagger a - [N]_{\alpha,\beta}), \label{genmacf}
\end{subeqnarray}
where $\dis [x]_{\alpha,\beta} = \frac{q^{\alpha x}- q^{\beta x}}
{q^\alpha - q^\beta}$ is a generalized $q$-bracket.
A similar generalization for the HY oscillator (GHY) gives
\begin{subeqnarray}
\lbrack a, a^{\dagger} \rbrack & = & [N + 1]_{\alpha,\beta} - 
[N]_{\alpha,\beta} \\
\ca_{4} & = & a^\dagger a - [N]_{\alpha,\beta}. \label{genhyo}
\end{subeqnarray}

We next consider realization of the $q$-deformed $su(2)$
algebra constructed from two independent
$q$-oscillators, $a, a^\dagger, N_a$ and $b, b^\dagger, N_b$.
Following ref \cite {mac,bied,kulish}
\begin{subeqnarray}
J_+  =  a^{\dagger}b, & & 
J_-  =  b^{\dagger}a, \\
J_0  =  \frac{1}{2} (N_a - N_b), & &
\ca =  \frac{1}{2} (N_a + N_b). \label{sch}
\end{subeqnarray}
Using the algebra defined in eq(\ref{macf}), we easily check that the
operators $J_\pm$ and $J_0$ obey the commutation relations:
\begin{subeqnarray}
{[}J_{\pm}, J_{0}{]} & =& \mp J_{\pm} \\
{[} J_{+}, J_{-} {]} & =&
\{- \ca_1 (q - q^{-1}) + 1 \} [2 J_0] 
\label{maccontr}
\end{subeqnarray}
Note that if we set $\ca_1=0$, we obtain the result
in ref \cite{mac,bied}. 
However, if we try to construct the realization using the
algebra defined in eq(\ref{hyo}), we arrive at the Fujikawa algebra
\cite{fuji2} with eq(\ref{maccontr}b) replaced by:
\begin{eqnarray}
{[} J_{+}, J_{-} {]} &=&
[2J_0] + \ca_2 \{ [ \ca - J_0 + 1] \nonumber \\
& & \hspace{-5mm} - [\ca - J_0] -
[\ca + J_0 + 1] + [\ca
+ J_0] \}.
\label{fujio}
\end{eqnarray}
This is not the conventional $q$-deformed $su(2)$ algebra as defined in
eq(\ref{cr}) unless $\ca_2 =
0$, which is the case in a Fock space representation.

Analogous Schwinger construction for the GMB 
and GHY algebras given by eq(\ref{genmacf}) and eq(\ref{genhyo})
can be constructed.
The commutation relations for the operators $\{ J_+, J_- \}$ for the
the GMB and GHY algebra are respectively
\be
{[} J_{+}, J_{-} {]} =
\{ \ca_3(q^{\alpha} - q^{\beta}) + 1 \} 
\frac{q^{\alpha N_a + \beta N_b} - 
q^{\alpha N_b + \beta N_a}}{q^\alpha - q^\beta}
\label{gencon1}
\ee
and
\begin{eqnarray}
{[} J_{+}, J_{-} {]} &= & 
\ca_4 \{ [N_{b} +1]_{\alpha,\beta} -[N_{b}]_{\alpha,\beta}     - [N_{a} + 1]_{\alpha,\beta} 
\nonumber \\
& & \mbox{\hspace{-5mm}} + [N_{a}]_{\alpha,\beta} \} +
  \frac{q^{\alpha N_a + \beta N_b} - q^{\alpha N_b + \beta N_a}}
{q^\alpha - q^\beta}.
\label{gencon2}
\end{eqnarray}
Note that when $\beta = -\alpha$, the term 
$\dis \frac{q^{\alpha N_a + \beta N_b} - q^{\alpha N_b + \beta N_a}}
{q^\alpha - q^\beta}
$ in eq(\ref{gencon1}) and eq(\ref{gencon2}) 
becomes $[2J_0]_{\alpha,\beta}$. 

One can also define 
operators $\Jp, \Jm $ and $\Jz$ using the relations 
\begin{subeqnarray}
\Jp  =  q^{-(\alpha + \beta)N_b} a^{\dagger}b, & & 
\Jm  =  b^{\dagger}a q^{-(\alpha + \beta)N_b}, \\
\Jz  =  \frac{1}{2} (N_a - N_b), & &
\ct =  \frac{1}{2} (N_a + N_b). \label{gensch}
\end{subeqnarray}
A straightforward calculation 
for the GMB oscillator
algebra and GHY oscillator yields 
\be
\Jp \Jm - q^{-(\alpha + \beta)} \Jm \Jp 
= \ca_3 \{ (q^{\alpha} - q^{\beta})  
+ 1 \} [ 2\Jz ]_{\alpha,\beta}.
\label{gencon3}
\ee
and
\begin{eqnarray}
& &  \Jp \Jm - q^{-(\alpha + \beta)} \Jm \Jp  \nonumber \\
& = &  q^{-(\alpha + \beta)N_b} \ca_4 \{ [ \ct - \Jz + 1] 
- [\ct - \Jz] -
[\ct + \Jz + 1] 
\nonumber \\
& & + [\ct
+ \Jz] \} 
+ [ 2\Jz]_{\alpha,\beta} 
\label{gencon4}
\end{eqnarray}
respectively.
In the above Schwinger construction, the expression
$[ 2\Jz ]_{\alpha,\beta}$ is obtained with a redefinition of the
commutation relation for the operators $\Jp$ and $\Jm$.
Note that for $\alpha + \beta = 0$, eqs(\ref{gencon3}) and
(\ref{gencon4}) reduce to eqs(\ref{gencon1}) and (\ref{gencon2}) respectively.

\section{$q$-Holstein-Primakoff Transformation} \label{qhols}
It is well-known that one can realize the undeformed $su(2)$ algebra
nonlinearly 
with one harmonic oscillator using the HP transformation.
A $q$-analogue of the transformation has also been studied\cite{chai}. 
The $q$-analogue of the HP transformation is defined by
the relations
\begin{subeqnarray}
J_+ & = & a^\dagger \sqrt{[2j - N]}, \\
J_- & = & \sqrt{[2j - N]} a, \\
J_0 & = & N - j, \label{hol}
\end{subeqnarray}
where $j$ is some $c$-number. 

It can be checked easily that under MB $q$-deformed
oscillators, the realization (\ref{hol}) leads to
\begin{subeqnarray}
\lbrack J_0 , J_\pm \rbrack & = & \pm J_\pm, \\
\lbrack J_+, J_- \rbrack & = & [2J_0] + \ca_1 q^{-2 J_0}; \label{hol1}
\end{subeqnarray}
whereas under HY oscillators, the commutation relations become
\begin{subeqnarray}
\lbrack J_0 , J_\pm \rbrack & = & \pm J_\pm, \\
\lbrack J_+, J_- \rbrack & = & [2J_0] + \ca_2 \{ [2j- N +1] 
\nonumber \\
& & \hspace{1cm}
- [2j - N] \} \\
& = & [2J_0] + \ca_2 \{ [j - J_0 + 1] 
\nonumber \\
& & \hspace{1cm}
- [j - J_0] \}. \label{hol2}
\end{subeqnarray}
It is interesting to compare eq(\ref{hol1}) and eq(\ref{hol2}) with
eq(\ref{cr}) and eq(\ref{fujio}) respectively. 

For the GMB and GHY oscillator
algebras defined by eqs(\ref{genmacf}) and (\ref{genhyo}),
one can also define the $q$-analogue of the HP
transformations in the most obvious manner by replacing the usual $q$-bracket
by its generalized $q$-bracket.  It turns out that 
the generalized $q$-HP 
transformations are then given by the relations
\begin{subeqnarray}
\Jp & = & q^{- \frac{\alpha + \beta}{2} N}
a^\dagger \sqrt{[2j - N]_{\alpha,\beta}}, \\
\Jm & = & \sqrt{[2j - N]_{\alpha,\beta}} a 
q^{- \frac{\alpha + \beta}{2} N}, \\
\Jz & = & N - j. \label{genhol}
\end{subeqnarray}

One easily verifies that under the GMB $q$-deformed 
oscillator, the realization turns out to be given by the relations
\begin{subeqnarray}
& & \lbrack \Jz , \Jpm \rbrack = \pm \Jpm,  \\
& & \Jp \Jm - q^{\alpha + \beta} \Jm \Jp  \nonumber \\
& = & [-2\Jz]_{\alpha,\beta} + \ca_3 q^{-2 \Jz \beta}; \label{genhol1}
\end{subeqnarray}
whereas under the GHY algebra, the same computation leads 
to the relations
\begin{subeqnarray}
& & \lbrack \Jz , \Jpm \rbrack = \pm \Jpm, \\
& & \Jp \Jm - q^{\alpha + \beta} \Jm \Jp  \nonumber \\
& = & q^{-(\alpha + \beta)N} \{ [2j - N + 1]_{\alpha,\beta} 
- [2j - N]_{\alpha,\beta} \} \ca_4 \nonumber \\ 
& & \hspace{30mm} + [-2 \Jz]_{\alpha,\beta}. \label{genhol2}
\end{subeqnarray}

\section{Contraction} \label{sect4}
So far we have tried to construct the $q$-deformed $su(2)$ from
$q$-oscillator algebras.  A somewhat reverse process, known as contraction,
is possible in
general.  For the undeformed case, we know that the
transformation\cite{gil} 
\be
\left(
\begin{array}{c}
h_+ \\
h_- \\
h_0 \\
1_h
\end{array}
\right) =
\left(
\begin{array}{cccc}
\mu & 0 & 0 & 0 \\
0 & \mu & 0 & 0 \\
0 & 0 & 1 & \frac{\eta}{2\mu^2}\\
0 & 0 & 0 & 1
\end{array}
\right)
\left(
\begin{array}{c}
J_+ \\
J_-\\
J_0 \\
\xi
\end{array}\right) \label{gilm}
\ee
maps the generators of $U(2)$, $J_\pm$ and $J_0$ with $[{\bf J}, \xi]=0$
under a change of basis to the generators $h_\pm, h_0$ and $1_h$ such that
\begin{subeqnarray}
\lbrack h_0, h_\pm \rbrack & = & \pm h_\pm \\
\lbrack h_+, h_- \rbrack & = & 2\mu^2 h_0 - \eta 1_h \\
\lbrack {\bf h}, 1_h] & = & 0. \label{contr}
\end{subeqnarray} 
One easily notes that the  commutation relations eq(\ref{contr})
are well-defined in the limit $\mu
\rightarrow 0$ despite the singularity in the
transformation. For $\mu \rightarrow 0$ and $\eta \rightarrow 1$,
the transformed algebra in eq(\ref{contr}) can be mapped isomorphically
to the standard oscillator algebra. This transformation is sometimes known as the
generalized In\"{o}n\"{u}-Wigner contraction.

The transformation given in eq(\ref{gilm}) allows for a simple extension
to the $q$-deformed case if one identifies the operators $\{ h_+, h_-,
h_0 \}$ as the operators $\{ a^\dagger, a, N^\prime \}$, the latter
satisfying the HY algebra with $\dis N^\prime = N + \frac{1}{2}$. 
Further one should also demand that
the operators $\{ J_+, J_-, J_0 \}$  obey the $q^{\frac{1}{2}}$-deformed
$su(2)$ algebra.  In particular, one can easily work out the
commutation relations for $[h_+,h_-]$ or equivalently $[a^\dagger, a]$
explicitly to get
\begin{eqnarray}
\lbrack h_+, h_- \rbrack & = & [a^\dagger, a] \nonumber \\
& = & \mu^2 [J_+, J_-] \nonumber \\
& = & \mu^2 \frac{q^{J_0} - q^{-{J_0}}}
{q^{\frac{1}{2}} -q^{-\frac{1}{2}}}  \nonumber \\
& = & \mu^2 \frac{q^{h_0}q^{-\frac{\eta}{2 \mu^2}\xi} -
q^{-h_0}q^{\frac{\eta}{2 \mu^2}\xi} }
{q^{\frac{1}{2}} -q^{-\frac{1}{2}}}.  \label{lhs}
\end{eqnarray}
However since  the
operators $\{ h, h^\dagger, h_0 \}$ or equivalently $\{ a, a^\dagger,
N^\prime \}$ obey the HY algebra, one can also work out the
commutation relation in eq(\ref{lhs}) in terms of the operator $h_0$. 
An straightforward computation yields
\begin{eqnarray}
\lbrack h_+, h_- \rbrack & = & [h_0 - \frac{1}{2}] -
[h_0 + \frac{1}{2}] \nonumber \\
& = & - \frac{q^{h_0} + q^{-h_0}}{q^{\frac{1}{2}} + q^{-\frac{1}{2}}}.
\label{rhs}
\end{eqnarray}
Consistency requirement for the  expressions in eq(\ref{lhs}) and
eq(\ref{rhs}) yields:
\begin{subeqnarray}
\frac{\mu^2}{q^{\frac{1}{2}} -q^{-\frac{1}{2}}} 
q^{-\frac{\eta}{2\mu^2}\xi} 
& = & - \frac{1}{q^{\frac{1}{2}} + q^{-\frac{1}{2}}}, \\
 \frac{\mu^2}{q^{\frac{1}{2}} -q^{-\frac{1}{2}}} q^{\frac{\eta}{2\mu^2}
\xi} & = &  \frac{1}{q^{\frac{1}{2}} + q^{-\frac{1}{2}}}. \label{coef}
\end{subeqnarray}
It is straightforward to solve
eq(\ref{coef}) for $\mu$ and $\eta\xi$ giving
\begin{subeqnarray}
\mu & = & \expon^{-i \frac{\alp}{2}}(\frac{q - 1}{q + 1})^{\frac{1}{2}} \\
\eta \xi & = & 2 \expon^{-i \alp} (\frac{q - 1}{q + 1}) \frac{i \alp}{\ln
q} \label{solve}
\end{subeqnarray}
where  $\alp = \frac{\pi}{2} +
\ell\pi$ ($\ell
\in {\bf Z}$) and we have appropriately chosen one branch when
taking the logarithm of complex number.

Thus, we observe that the relation obtained in ref \cite{cel1,oh2}
between HY oscillator and
$su_{\sqrt{q}}(2)$ algebra can be regarded as the $q$-analogue of the
transformation given in eq(\ref{gilm}) if we write
\be
\left(
\begin{array}{c}
a_+ \\
a \\
N^\prime \\
1
\end{array}
\right) =
\left(
\begin{array}{cccc}
\expon^{-i\frac{\alp}{2}} (\frac{q - 1}{q + 1})^{\frac{1}{2}}
& 0 & 0 & 0 \\
0 &
\expon^{-i\frac{\alp}{2}} (\frac{q - 1}{q + 1})^{\frac{1}{2}}
& 0 & 0 \\
0 & 0 & 1 & \frac{i}{\ln q}\\
0 & 0 & 0 & 1
\end{array}
\right)
\left(
\begin{array}{c}
J_+ \\
J_-\\
J_0 \\
\alp  1
\end{array}\right) \label{cele}
\ee
in which one easily identifies the quantities $\mu$, $\eta$ and $\xi$ in
eq(\ref{gilm}) by
$\dis \mu = \expon^{-i\frac{\alp}{2}} (\frac{q - 1}{q +
1})^{\frac{1}{2}} $, $\dis
\eta = \frac{2i \expon^{-i \alp}}{\ln q} (\frac{q - 1}{q + 1})$
and $\xi = \alp $. We would like to 
emphasize again that the operators $J_\pm$
and $J_0$ in this case
obey the $q^\frac{1}{2}$-deformed commutation relations
in eq(\ref{cr})\cite{cel1,oh2}.
In the
limit $q \rightarrow 1$, this transformation is again singular but again
the commutation relations for the oscillator algebra are well-defined and
become the undeformed oscillator algebra. Furthermore, for generic $q$,
the coproduct,
counit and antipodes for the $q$-deformed $su(2)$ carry directly through
the transformation, endowing the HY oscillator with a Hopf
structure. This Hopf structure however breaks down in the limit when $q
\rightarrow 1$ whereas the Hopf structure of $su_{\sqrt{2}}(2)$ becomes
cocommutative in the same limit. From refs \cite{oh2,fuji1}, it is not
difficult to show that the positive norm requirement and the truncation
condition for the states of the HY $q$-oscillator are in conflict
with each other.  Thus the HY $q$-oscillator algebra is not the
same as the $su_{\sqrt{q}}(2) $ algebra.

\begin{comment}
In particular, when $q= \expon^{i\epsilon}$, one finds that the
contraction in eq(\ref{cele}) takes the form:
\be
\left(
\begin{array}{c}
a_+ \\
a \\
N \\
1
\end{array}
\right) =
\left(
\begin{array}{cccc}
\expon^{-i\frac{\alp}{2}} \sqrt{\tan(\frac{\epsilon}{2})}
& 0 & 0 & 0 \\
0 &
\expon^{-i\frac{\alp}{2}} \sqrt{(\tan(\frac{\epsilon}{2})}
& 0 & 0 \\
0 & 0 & 1 & \frac{1}{\epsilon}\\
0 & 0 & 0 & 1
\end{array}
\right)
\left(
\begin{array}{c}
J_+ \\
J_-\\
J_0 \\
\alp  1
\end{array}\right). \label{cele2}
\ee
It is not difficult to check that the operators $\{ a, a^\dagger, N \}$
satisfy the HY oscillators in eq(\ref{hyo}) but not the
MB oscillator relations in eq(\ref{macf}).  One also
notes that since $\ln \exp^{i \epsilon}
\end{comment}

The MB oscillator algebra can be shown via the map $a = q^{N/2}
A$, $a^\dagger = A^\dagger q^{N/2}$ to be equivalent to the algebra
${\cal A}_q$ with
operators $\{ A, A^\dagger, N \}$ satisfying
\begin{subeqnarray}
\lbrack A, A^\dagger \rbrack & = & q^{-2N} \\
\lbrack N, A \rbrack & = & - A \\
\lbrack N, A^\dagger \rbrack & = & A^\dagger. 
\end{subeqnarray}
In fact, Chaichian and Kulish \cite{chai} have shown that the map
\begin{subeqnarray}
A & = & \lim_{s \rightarrow \infty} \frac{(q - q^{-1}) }{q^s} J_+ \\
A^\dagger & = & \lim_{s \rightarrow \infty} \frac{(q - q^{-1}) }{q^s} J_- \\
N & = & s - J_0 
\end{subeqnarray}
allows the contraction of $su_q(2)$ to the MB $q$-oscillator algebra.
Note that this contraction clearly lifts the
highest weight representation to infinity so that there exists an
infinite tower of states needed for the oscillator algebra ${\cal A}_q$.
Although this
contraction does not induce a coproduct for $\{ A, A^\dagger, N \}$,
it admits a coaction $\Psi:  {\cal A}_q
\rightarrow  {\cal A}_q \otimes SU_q(2)$ given by
\begin{subeqnarray}
\Psi(N) & = & N - J_0, \\
\Psi(A) & = & A q^{-J_0} + \sqrt{q - q^{-1}} q^{-N} J_+ ,\\
\Psi(A^\dagger) & = & A^\dagger q^{-J_0} + \sqrt{q - q^{-1}} q^{-N} J_-. 
\end{subeqnarray}
This coaction satisfies the associative axioms namely
\begin{subeqnarray}
(\Psi \otimes 1) \circ \Psi &=& (1 \otimes \Psi) \circ \Psi \\
(1 \otimes \epsilon) \circ \Psi & = & 1
\end{subeqnarray}
where $\epsilon$ is the counit.  Further, one easily checks that the
homomorphism axiom is consistent, namely
\be
\Psi([x, y])  =  [\Psi(x), \Psi(y)]
\ee
where $x, y \in \{ A, A^\dagger, N \}$. In the framework of In\"{o}ue-Wigner
transformation, there seems to be a singular transformation
\be
\left(
\begin{array}{c}
A \\
A^\dagger \\
N \\
1
\end{array}
\right) =
\left(
\begin{array}{cccc}
\dis \frac{\sqrt{q - q^{-1}}}{q^s} & 0 & 0 & 0 \\
0 & \dis \frac{\sqrt{q - q^{-1}}}{q^s} & 0 & 0 \\
0 & 0 & -1 & s\\
0 & 0 & 0 & 1
\end{array}
\right)
\left(
\begin{array}{c}
J_- \\
J_+\\
J_0 \\
1
\end{array}\right) .
\ee
The contraction from $su_q(2)$ to MB oscillator algebra occurs
in the singular limit $s \rightarrow \infty$, but in this case, the
natural coproduct for $su_q(2)$ does not survive in this limit.
This contraction is essentially similar to the one proposed by J. Ng \cite{ng}.
%%%%%%%%%%%%%%%%%%%%%%%%%%%%
A different contraction proposed by Celeghini et al \cite{cel1,cel2}
involves the transformation
\be
\left(
\begin{array}{c}
B \\
B^\dagger \\
N \\
H \\
\omega
\end{array}
\right) =
\left(
\begin{array}{ccccc}
\eta & 0 & 0 &  0  & 0 \\
0 & \eta & 0 & 0  & 0 \\
0 & 0 & -1 & \eta^{-2} & 0 \\
0 & 0 & 0 & 2 & 0 \\
0 & 0 & 0 & 0 & \eta^{-2}
\end{array}
\right)
\left(
\begin{array}{c}
J_+ \\
J_-\\
J_0 \\
K \\
\log q
\end{array}\right) .
\ee
where $K$ is the so-called $U(1)$ generator.
Under this transformation, the operators $\{ B, B^\dagger, N, H \}$
obey in the limit $\eta \rightarrow 0 $ the relations
\begin{subeqnarray}
\lbrack B, B^\dagger \rbrack & = & \dis \frac{\sinh (\frac{\omega H}{2})}{\frac{\omega}{2}}  \\
\lbrack N, B \rbrack =  - B, & \lbrack N, B^\dagger \rbrack
 =   B^\dagger, &  \lbrack H, N \rbrack  =  0 \\
\lbrack H, B \rbrack = & \lbrack H, B^\dagger \rbrack = & 0 
\end{subeqnarray}
This contraction induces a coalgebraic structure inherited from the
original Hopf algebra of $su_q(2)$.
\begin{comment}
, namely
\begin{subeqnarray}
\Delta(B) & = & \expon^{-\omega H /4} \otimes B + B \otimes \expon^{-\omega H
/4} \\
\Delta(B^\dagger) & = & \expon^{-\omega H /4} \otimes B^\dagger + B^\dagger
\otimes \expon^{-\omega H
/4} \\
\Delta (N) & = & 1 \otimes N + N \otimes 1 \\
\Delta (H) & = & 1 \otimes H + H \otimes 1.
\end{subeqnarray}
\end{comment}
The algebra generated by the operators $\{ B, B^\dagger, H, N \}$ is not quite
the $q$-deformed oscillator algebra although we can get the usual
undeformed oscillator in the limit $\omega \rightarrow 0$.

\section{Representations}

We can gain some insights into the the linear transformation which we
have encountered
in the previous section by looking more closely at a representation of
the HY oscillator algebra. To obtain a representation of the HY
algebra\cite{oh2}, we note
that $N$ commutes with $a^\dagger a$ and $a a^\dagger$.  As a result we
can construct a vector $|\psi_0>$ which is a simultaneous eigenstate of
$N$ and $a^\dagger a$ so that
\begin{subeqnarray}
N |\psi_0> & = & \nu_0 |\psi_0> \\
a^\dagger a |\psi_0> & = & \lambda_0 |\psi_0>
\end{subeqnarray}
where $\nu_0$ and $\lambda_0$ are the corresponding eigenvalues. We
shall further assume that the operator $N$ is Hermitian so that its
eigenvalue $\nu_0$ is real.
\begin{comment}
If we
demand hermiticity for the operator $N$ then one can show that $|q| =1$.
This is not necessary the case if we drop the hermiticity condition.
Here, we shall not demand for
hermiticity of $N$ and allow $\nu_0$ to be arbitrary.
\end{comment}

From the eigenstate, $|\psi_0>$, one can construct other eigenstates of
$N$ by defining
\begin{subeqnarray}
|\psi_n> & = & (a^\dagger)^n |\psi_0> \\
|\psi_{-n}> & = & a^n |\psi_0> 
\end{subeqnarray}
for some positive integer $n$. With these definitions, one easily shows
that 
\begin{subeqnarray}
a^\dagger |\psi_n> & = & |\psi_{n+1}> \\
a^\dagger |\psi_{-n}> & = & \lambda_{-n+1}|\psi_{-n+1}> \\
a |\psi_n> & = & \lambda_n |\psi_{n-1}> \\
a |\psi_{-n}> & = & |\psi_{-n-1}> \\
N |\psi_{\pm n}> & = & (\nu_0 \pm n) |\psi_{\pm n}> 
\end{subeqnarray}
where 
\begin{subeqnarray} 
\lambda_n & = & \lambda_0 + \frac{q^{\frac{1}{2}n} -
q^{-\frac{1}{2}n}}{q^{\frac{1}{2}} - q^{-\frac{1}{2}}}  \frac{q^{\nu_0 +
\frac{n}{2}} + q^{-\nu_0 - \frac{n}{2}}}{q^{\frac{1}{2}} +
q^{-\frac{1}{2}}} \\
& =& \lambda_0 + [n + \nu_0] - [\nu_0]. \label{lambdan}
\end{subeqnarray}
Note that the oscillator algebra still admits an infinite number of
states and the representation at this stage is different from
$\U_q(su(2))$ whose finite-dimensional
representation requires a highest weight state. 
One then imposes a truncation on the tower of states and set $a |\psi_0> =
0$ giving $\lambda_0 = 0$ and $|\psi_{-n}>= 0$ for any $n > 0$. Let
$|\psi_k>$ be the highest weight state so that 
$a^\dagger |\psi_k> = 0$ with integer $k > 0$. 
Since $\ca_2=a^\dagger a - [N] = a a^\dagger - [N + 1]$, one finds by
considering the action of
$\ca_2$ on $|\psi_k>$ that the following condition must be satisfied:
\be
[\nu_0 + k + 1] = [\nu_0] .\label{trun}
\ee
For real $q$, $k = -1$ is the only solution, but this is not acceptable.
However, for complex
$q$ with  $|q| = 1$, truncation is possible.
It is not difficult to solve eq(\ref{trun}) for $\nu_0$ in this case. 
Writing $q = \expon^{i \epsilon}$, 
one can show that for arbitrary $\epsilon$, eq(\ref{trun}) leads to
\be
\nu_0 \epsilon = \frac{-(k+1)\epsilon}{2} + (\ell + \frac{1}{2}) \pi,
\mbox{\hspace{3cm}} \ell \in {\bf Z}
\label{soln}
\ee
%We note here that the 
This
result needs not be consistent with the
condition for positivity of norms \cite{oh2,fuji1} which by eq(\ref{lambdan})
is
\be
\lbrack n + \nu_0  \rbrack - \lbrack \nu_0 \rbrack \geq 0 \label{positive}
\ee for all integers $n \leq k $.
To see this, we substitute eq(\ref{soln}) into
the left hand side of condition (\ref{positive}) and see that
$$
\lbrack n + \nu_0 \rbrack - [\nu_0] = \frac{(-1)^\ell}{\sin \epsilon}
\{  \cos(\frac{k + 1}{2} - n) \epsilon - \cos \frac{k + 1}{2} \epsilon \}
$$
which needs not be positive for arbitrary  $\epsilon$.
This means that for arbitrary $\epsilon$, we cannot proceed to identify
the HY oscillator algebra with $su_{\sqrt{q}}(2)$ algebra.  To
identify the two algebras, we have to truncate the tower of states of
the HY oscillator algebra. However, truncation and positive norm
requirement can both be satisfied only for certain value of $\epsilon$. 
In short, the HY oscillator algebra and $su_{\sqrt{q}}(2)$ algebra are
equivalent only for certain $q$-values.

\vspace{1cm}

\centerline{\large \bf Acknowledgments} 

\noindent We wish to thank Prof. Kazuo Fujikawa for many helpful 
suggestions and discussions.

\end{document}